\documentclass[fleqn]{article}
\usepackage{amsmath}
\usepackage{amssymb}

\newcommand{\openR}{\mathbb{R}}
\newtheorem{theorem}{Theorem}[section]
\newtheorem{lemma}[theorem]{Lemma}
\newtheorem{proposition}[theorem]{Proposition}
\newtheorem{corollary}[theorem]{Corollary}
\newtheorem{remark}[theorem]{Remark}

\newcommand{\cov}{\mbox{Cov}\,}

\newcommand{\E}{\mathbb{E}\,}
\newcommand{\eps}{\varepsilon}
\newcommand{\bdots}{\begin{array}{c}
         \vspace*{2pt}\hspace*{8pt} \cdot \\ \vspace*{-22pt} \\
         \vspace*{2pt}\hspace*{4pt} \cdot \hspace*{5pt} \\ \vspace*{-22pt} \\
         \vspace*{3.5pt}\cdot \hspace*{10pt} \end{array}}

\def\dddots{\mathinner{\mkern1mu\raise1pt\vbox{\kern7pt\hbox{.}}
\mkern2mu\raise4pt\hbox{.}\mkern2mu\raise7pt\hbox{.}\mkern1mu}}
 \def\d4dots{\mathinner{\vbox{\kern7pt\hbox{.}}
    \mkern2mu\raise3pt\hbox{.}
    \mkern2mu\raise6pt\hbox{.}
    \mkern2mu\raise9pt\hbox{.}}}

\begin{document}
\pagestyle{empty}
 \Large
\begin{center}
The Bezoutian and Fisher's information matrix \\ of an ARMA
process
\bigskip \\
Andr\'{e} Klein\footnote{Department of Quantitative Economics,
 Universiteit van Amsterdam,
Roetersstraat 11, 1018 WB Amsterdam, {\tt a.a.b.klein@uva.nl}} and
Peter Spreij\footnote{ Korteweg-de Vries Institute for
Mathematics, Universiteit van Amsterdam, Plantage Muidergracht 24,
1018 TV Amsterdam, {\tt spreij@science.uva.nl}}
\end{center}

\normalsize
\begin{abstract}
\noindent In this paper we derive some properties of the Bezout
matrix and relate the Fisher information matrix for a stationary
ARMA process to the Bezoutian. Some properties are explained via
realizations in state space form of the derivatives of the white
noise process with respect to the parameters. A factorization of
the Fisher information matrix as a product in factors which
involve the Bezout matrix of the associated AR and MA polynomials
is derived. From this factorization we can characterize
singularity of the Fisher information matrix.

\begin{description}
\item{\sl Keywords:} ARMA process, Fisher information matrix, Stein's equation,
Syl\-ves\-ter's resultant matrix, Bezout matrix, state space
realization
\item{\sl AMS subject classification:} 15A23, 15A24, 15A57,  62M10,
93B15, 93B17
\end{description}

\end{abstract}

\normalsize

\newpage
~
\newpage
\pagestyle{plain} \setcounter{page}{1}

\section{Introduction}
The Cram\'{e}r-Rao lower bound on the covariance matrix of an
estimator is a classical result in statistics, see Cram\'{e}r
\cite{cramer} and Rao \cite{rao}. This bound is given by the
inverse of Fisher's information matrix. For regular statistical
models, it is also known that the maximum likelihood estimator is
asymptotically normal with this inverse as the asymptotic
covariance matrix.  Therefore it is natural to ask for conditions
of an underlying statistical model that guarantee non-singularity
of this matrix. In the present paper we are concerned with the
Fisher information matrix for (stationary) autoregressive moving
average (ARMA) models.  The information matrix is singular in the
presence of common roots of the AR and the MA polynomial and vice
versa. This fact is considered to be well-known in time series
analysis, see \cite{mcleod} or \cite {potscher} for an early
discussion of this phenomenon, and~\cite{kleinspreij} for the
extension to ARMA models with an exogenous input (the ARMAX case).

In~\cite{gyor} properties of the Fisher information matrix for an
ARMA process have been derived using contour integration in the
complex plane and state space realizations of the ARMA process
itself. In the present paper we study Fisher's information matrix
by means of state space realizations for the score process and by
linking Fisher's information matrix to the Sylvester resultant
matrix and the Bezout matrix associated with the autoregressive
and moving average polynomials.

The role of the resultant matrix has been discussed in various
studies in the fields of time series and systems theory. For
instance, in \cite {astromsoderstrom} this matrix shows up in a
convergence analysis of maximum likelihood estimators of the ARMA
parameters (more precisely in the study of the convergence of the
criterion function), in Barnett \cite{barnett} a relationship
between Sylvester's resultant matrix and the companion matrix of a
polynomial is given. Kalman \cite{kalman} has investigated the
concept of observability and controllability in terms of
Sylvester's resultant matrix. Similar results can be found in
Barnett \cite{barnett2}, which contains further discussions and
references on these topics. But, it seems that the use of the
Bezout matrix has not been recognized yet. For ARMA models we will
show that Fisher's information matrix can be factorized, where one
of the factors is expressed in terms of the Bezout matrix. Also
from this it follows that Fisher's information matrix is singular
if and only if the AR and MA polynomials have a non-trivial
greatest common divisor. Singularity of the information matrix can
thus be interpreted as the result of overparametrization of the
chosen ARMA model and of using a model of too high order. In
S\"{o}derstr\"{o}m \& Stoica \cite[pages 162 ff.]{soderstrom} a
discussion on over\-pa\-rametri\-zation in terms of the transfer
function of a system can be found. In a static context, Fisher's
information matrix has already been studied in \cite{rothenberg}
for problems of local and global identifiability.

The paper is organized as follows. In section~\ref{sec:state} the
main results are state space realizations for the derivatives of
the noise process and properties of these realizations are
presented. In section~\ref{sec:bezout} we study some properties of
the Bezout matrix as well as its kernel.
Section~\ref{sec:armabezout} is devoted to further properties of
the Bezout matrix, to be exploited in sections~\ref{sec:stein}
and~\ref{section:fisher}. In the first of these sections, we study
singularity of solutions to certain Stein equations with
coefficients related to the AR and MA polynomials, whereas in
section~\ref{section:fisher} all previous results are assembled to
characterize non-singularity of Fisher's information matrix.

\section{Computations in state space}\label{sec:state}

\setcounter{equation}{0}

Consider the
following two scalar monic polynomials in the variable $z$.
\begin{eqnarray}
\hat{a}(z) & = & z^p+a_1z^{p-1}+\cdots+a_p \nonumber \\
\hat{c}(z) & = & z^q+c_1z^{q-1}+\cdots+c_q. \nonumber
\end{eqnarray}
By $a$ and $c$ we denote the reciprocal polynomials, so
$a(z)=z^p\hat{a}(z^{- 1})$ and $a(z)=z^qc(z^{-1})$, and also write
$a^\top =(a_1,\ldots,a_p)$ and $c^\top =(c_1,\ldots,c_q)$. Usually
no confusion between the notation $a$ for the polynomial and
vector will arise, but sometimes we will write $a(\cdot)$ when a
polynomial is considered.
\medskip\\
Consider  the stationary
ARMA($p$,$q$) process $y$ that satisfies
\begin{equation}
a(L)y=c(L)\varepsilon  \label{eq:arma}
\end{equation}
with $L$ the lag operator and $\varepsilon$ a white noise
sequence. We make the assumption (to give the expressions that we
use below the correct meaning) that $y$ is {\it causal} and {\it
invertible}, i.e.\ both ${a}$ and ${c}$ have
only zeros outside the unit circle (equivalently, $\hat{a}$ and $\hat{c}$ have only zeros inside the unit circle). \\
Let $\theta=(a_1,\ldots,a_p,c_1,\ldots,c_q)$ and denote by
$\varepsilon^{\theta_i}_t$ the derivative of $\varepsilon_t$ with
respect to $\theta_i$. Then we obtain by differentiation of
(\ref{eq:arma}) the formal expressions
\begin{eqnarray}
\varepsilon^{a_j}_t & = & \frac{1}{a(L)}\varepsilon_{t-j} \label{eq:epsa} \\
\varepsilon^{c_l}_t & = & -\frac{1}{c(L)}\varepsilon_{t-l}.
\label{eq:epsc}
\end{eqnarray}
Let $\dot{\varepsilon}_t=\dot{\varepsilon}_t(\theta)$ denote the
row vector with elements $\varepsilon^{\theta_i}_t$ and
$\xi_t=\dot{\varepsilon}_t^\top$. See section~\ref{section:fisher}
for the relation with the stationary Fisher information matrix of
the ARMA process $y$. We introduce some auxiliary notation. Write
for each positive integer $k$
\begin{eqnarray*}
u_k(z) & = & (1,z,\ldots,z^{k-1})^\top  \\
u^*_k(z) & = & (z^{k-1},\ldots,1)^\top =z^{k- 1}u_k(z^{-1}).
\end{eqnarray*}
Let us compute the transfer function $\tau(z)$ that relates $\xi$
to $\varepsilon$ by replacing $L$ with $z^{-1}$ in
equations~(\ref{eq:epsa}) and~(\ref{eq:epsc}). Here $z^{-1}$
represents the forward shift. One obtains from~(\ref{eq:epsa})
and~(\ref{eq:epsc})
\begin{equation}\label{eq:tau}
\tau(z)=
\begin{pmatrix}
\frac{1}{\hat{a}(z)}u^*_p(z) \\
-\frac{1}{\hat{c}(z)}u^*_q(z)
\end{pmatrix}.
\end{equation}
In~\cite{gyor} we have investigated certain controllable or
observable realizations of the ARMA process $y$. There we have
also briefly outlined a procedure without detailed proofs to
obtain from these realizations also realizations for the process
$\dot{\varepsilon}$. We repeat the conclusions, but give in the
present paper a short proof of them, based on transfer function
considerations, without using the realizations of the ARMA process
$y$ itself. Let $e_n$ be the first basis vector of the Euclidean
space $\mathbb{R}^n$. When no confusion can arise (often in the
proofs), we often simple write $e$, which we will also use as the
notation for the first basis vector in Euclidean spaces of
different dimensions.  By $J$ we denote the forward shift matrix,
$J_{ij}=1$ if $i=j+1$ and zero else.  \\
Similarly, we denote by $I$ the identity matrix of the appropriate
size and $0$ stands for the zero vector or matrix of appropriate
dimensions. Occasionally these matrices and vectors will have a
subscript, when it is
necessary to indicate the sizes.  \\
Let
$\hat{g}=\hat{a}(z)\hat{c}(z)=z^{p+q}+\sum_{i=1}^{p+q}g_iz^{p+q-
i}$ and  $\hat{g}$  the vector
$\hat{g}=(g_{p+q},\ldots,g_{1})^\top$. Likewise we write
$g(z)=a(z)c(z)$ and $g=(g_1,\ldots,g_{p+q})$.
\medskip\\
The Sylvester resultant matrix $R$ of $\hat{c}$ and $-\hat{a}$ is
defined as the $(p+q)\times(p+q)$ matrix
\begin{equation}\nonumber
R(c,-a) =
  \begin{pmatrix}
    R_p(c) \\
    -R_q(a)
  \end{pmatrix},
\end{equation}
where $R_p(c)$ is the $q\times(p+q)$ matrix
\begin{equation}\nonumber
R_p(c)=
  \begin{pmatrix}
  1      & c_1    & \cdots &  c_q   &  0   & \cdots &  0   \\
  0      & 1      & c_1    & \cdots &  c_q &        &      \\
  \vdots & \ddots & \ddots & \ddots &      & \ddots &      \\
  0      & \cdots & 0      & 1      & c_1  & \cdots &  c_q
  \end{pmatrix}
\end{equation}
and $R_q(a)$ is the $q\times(p+q)$ matrix given by
\begin{equation}\nonumber
R_q(a) =
  \begin{pmatrix}
  1      & a_1    & \cdots &  a_p   &  0   & \cdots &  0   \\
  0      & 1      & a_1    & \cdots &  a_p &        &      \\
  \vdots & \ddots & \ddots & \ddots &      & \ddots &      \\
  0      & \cdots & 0      & 1      & a_1  & \cdots &  a_p
  \end{pmatrix}.
\end{equation}
In the presence of common roots of $\hat{a}$ and $\hat{c}$ the
matrix $R(c,-a)$ becomes singular. Moreover it is known (see e.g.
\cite[page 106] {vanderwaerden}) that
\begin{equation}
\det R(c,-a) = (-1)^p\prod_{i=1}^p \prod_{j=1}^q (\gamma_j -
\alpha_i) \label{eq:det}
\end{equation}
where the $\alpha_i$ and the $\gamma_j$ are the roots of $\hat{a}$
and $\hat{c}$ respectively.
\medskip\\
Next we introduce the matrices $F$ and $G$ defined by
\begin{equation}\label{eq:F}
F=\begin{pmatrix}
J-e_pa^\top & 0 \\
0 & J-e_{q}c^\top
\end{pmatrix}
\end{equation}
and
\begin{equation}\label{eq:G}
G= J-e_{p+q}g^\top.
\end{equation}
\begin{lemma}\label{lemma:RR}
Let $F$ and $G$ be as in~(\ref{eq:F}) and (\ref{eq:G}). Then the
following relation holds.
\begin{equation}
R(c,-a)G= F R(c,-a). \label{eq:RR}
\end{equation}
\end{lemma}
{\bf Proof.} The easiest way to see this, is to multiply both
sides of this equation on the right with $u_{p+q}(z)$. Then we
compute on the left hand side the product
\begin{align*}
R(c,-a)(J-e_{p+q}g^\top)u_{p+q}(z) & =   R(c,-a)(u_{p+q}(z)-g(z)e) \\
& =
\begin{pmatrix}
c(z)u_p(z)-g(z)e_p \\
-a(z)u_q(z)+g(z)e_q
\end{pmatrix}.
\end{align*}
The computations on the right hand side are of a similar nature
and an easy comparison yields the result.
\hfill$\square$\medskip\\
We now present the first realization of the process $\dot{\varepsilon}$.
\begin{proposition}\label{prop:contrel}
The process $\xi=\dot{\varepsilon}^\top$ can be realized by the
following stable and controllable system
\begin{align}
Z_{t+1} & =GZ_t+e\varepsilon _t  \label{eq:obsreal} \\
\xi_t & =CZ_t,  \label{eq:obsreall}
\end{align}
where $G$ is as in~(\ref{eq:G}) and $C=R(c,-a)$. This system is
observable iff the polynomials $a$ and $c$ have no common zeros.
\end{proposition}
{\bf Proof.} Let us compute the transfer function $\tau$ of the
above system. Standard computations show that $(z-G)^{-1}
e=\frac{1}{\hat{g}(z)}u^*_{p+q}(z)$. The trivial observations
$R_p(c)u^*_{p+q}(z)=\hat{c}(z)u^*_p(z)$ and
$R_q(a)u^*_{p+q}(z)=\hat{a}(z)u^*_q(z)$ then lead to the
conclusion that $\tau(z)=C(z-G)^{-1}e$ is exactly the same as
in~(\ref{eq:tau}). The system is obviously controllable. The
observability matrix of the system involves products of the form
$R(c,-a)G^k$ ($k=0,\ldots,p+q-1$). In view of lemma~\ref{lemma:RR}
these can be written as $F^kR(c,-a)$, from which the assertion
follows. Stability is an immediate consequence of the assumptions
on the polynomials $a$ and $c$. Indeed, the characteristic
polynomial of $G$ is $\hat{g}=\hat{a}\hat{c}$, which has its zeros
inside the unit circle.
\hfill$\square$\medskip\\
An alternative (observable) realization of the process
$\dot{\varepsilon}$ is given in the next proposition.
\begin{proposition}\label{prop:eps2}
The process $\xi=\dot{\varepsilon}^\top$ is the state process of
the stable system given by
\begin{equation}
\xi_{t+1}=F\xi_t +B\varepsilon _t, \label{eq:epsstate}
\end{equation}
where $F$ is as in~(\ref{eq:F}) and $B=\begin{pmatrix}
e_p \\
-e_q \end{pmatrix}$. This system is controllable iff $\hat{a}$ and
$\hat{c}$ have no common zeros.
\end{proposition}
{\bf Proof.} Again the proof that this realization produces
$\dot{\varepsilon}$ boils down to computing the transfer function,
like we did in the proof of proposition~\ref{prop:contrel}. The
computations needed for this have been encountered there, so we
skip them. To explain the statement on controllability, we
consider the equation ($u,v$ are row vectors and $\lambda$ is an
arbitrary complex number)
\[
(u,v)\begin{pmatrix}
F_a-\lambda & 0 & e\\
0 & F_c-\lambda & -e \end{pmatrix}=0,
\]
where $F_a=J-e_pa^\top$ and $F_c=J-e_qc^\top$. This equation is
equivalent to $u(F_a-\lambda )=0$, $ v(F_c-\lambda)=0$ and
$(u-v)e=0$. We first consider the case where $\hat{a}$ and
$\hat{c}$ have no common zeros. Suppose that $u=0$. Then we have
that $v(F_c-\lambda)=0$ and $ve=0$. Since $(F_c,e)$ is
controllable, $v$ must be zero as well. Therefore we will assume
that there is a nonzero solution $u$. Then $\lambda$ must be a
root of $\hat{a}(z)=0$. If $v=0$, then we also have $ue=0$. This
situation cannot happen since $(F_a,e)$ is a controllable pair. So
we have to assume that $v$ is not zero, but then $\lambda$ is also
a root of $\hat{c}(z)=0$.  It then follows from the above that
this cannot happen. Hence $(F,B)$ is controllable. In the other
case $\hat{a}$ and $\hat{c}$ have a common zero $\lambda$. In this
case $u$ is the row vector
$(1,\hat{a}_1(\lambda),\ldots,\hat{a}_{p-1}(\lambda))$, where the
$\hat{a}_i$ are the H\"orner polynomials, defined by
$\hat{a}_0(z)=1$, $\hat{a}_k(z)=z\hat{a}_{k-1}(z)+a_k$ and we have
a similar expression for $v$. One obviously then also has
$(u-v)e=0$ and hence the system is not controllable. Stability
follows upon noting that the characteristic polynomial of $F$ is
equal to $\hat{g}$. \hfill$\square$
\begin{remark}
{\em By lemma~\ref{lemma:RR}, the realization of
proposition~\ref{prop:eps2} is connected to the one in
proposition~\ref{prop:contrel} in a very simple way. Starting from
equation~(\ref{eq:obsreall}), one obtains
\begin{align*}
\xi_{t+1}   & = R(c,-a)Z_{t+1} \\
  & = R(c,-a)(GZ_t+e\eps_t) \\
  & =  F\xi_t + B\eps_t.
\end{align*}
 }
\end{remark}
\begin{remark}\label{rem:independent}
{\em Notice that the realizations of propositions~\ref{prop:eps2}
and~\ref{prop:contrel} illustrate the well known fact that
$\dot{\eps}_t$ depends on $\eps_s$ for $s<t$ only, and hence is
stochastically independent of $\eps_t$.
}
\end{remark}

\section{The Bezoutian}\label{sec:bezout}

\setcounter{equation}{0}

We follow the notation of Lancaster \&
Tismenetsky~\cite{lancastertismenetsky}. Recall the following
definitions. In this section and henceforth we assume that $p$ and
$q$ are taken to have a common value, denoted by $n$, to yield
many of the subsequent expression meaningful. We consider
polynomials $a$ and $b$ given by $a(z)=\sum_{k=0}^na_kz^k$ and
$b(z)=\sum_{k=0}^nb_kz^k$. We will always assume that the constant
term $a_0=1$ and likewise for $b$ and other polynomials.
\medskip\\
Consider the Bezout matrix $B(a,b)$ of the polynomials $a$ and
$b$. It is defined by the relation
\begin{equation}\nonumber
a(z)b(w)-a(w)b(z) = (z-w)u_n(z)^\top B(a,b)u_n(w).
\end{equation}
We also often call this matrix the Bezoutian. Introduce for a
given complex number $\phi$ the matrices $U_\phi$ as follows.
\[ U_\phi =
\begin{pmatrix}
1      & 0      & \ldots & \ldots & 0 \\
-\phi  & 1      &        &        & \vdots \\
0      &        & \ddots &        & \vdots \\
\vdots &        &        & \ddots & 0 \\
0      & \ldots & 0      & -\phi  & 1
\end{pmatrix}.
\]
We also need the inverses $T_\phi$ of the matrices $U_\phi$. These
take the form
\[ T_\phi=
\begin{pmatrix}
1      & 0      & \ldots & \ldots & 0 \\
\phi   & 1      & 0      &        & \vdots \\
\phi^2 & \phi   & \ddots & \ddots & \vdots \\
\vdots &        & \ddots & \ddots & 0 \\
\phi^{n-1} & \ldots & \phi^2 & \phi  & 1
\end{pmatrix}.
\]
Observe that matrices $U_\phi$ and  $U_\psi$ commute, as well as
$T_\phi$ and $T_\psi$.
\medskip\\
Consider again $a$ and $b$, $n$-th order polynomials with constant
term equal to 1. Let $(1-\alpha_1z)$ be a factor of $a(z)$ and
$(1-\beta_1z)$ be a factor of $b(z)$. Of course, $\alpha_1$ and
$\beta_1$ are zeros of $\hat{a}$ and $\hat{c}$ respectively. Write
$a(z) = (1-\alpha_1z)a_{-1}(z)$ and $b(z)= (1-\beta_1z)b_{-
1}(z)$. Continuing this way, for $\alpha_1,\ldots,\alpha_n$ we
define recursively $a_{-(k-1)}(z)=(1-\alpha_kz)a_{-k}(z)$ and
polynomials $b_{-k}$ similarly. We also put $a_0(z)=a(z)$ and
$b_0(z)=b(z)$. The following proposition is not completely
necessary for what follows, but may be of independent interest.
\begin{proposition}\label{prop:bezexp}
With the above introduced notation we have
\begin{align} \frac{a(z)b(w)-a(w)b(z)}{z-w}  = &
(1-\alpha_1z)(1-\beta_1w)\frac{a_{-1}(z)b_{-
1}(w)-a_{-1}(w)b_{-1}(z)}{z-w}\nonumber \\
& \mbox{}+ (\beta_1-\alpha_1)a_{-1}(w)b_{-1}(z).\label{eq:bezdec1}
\end{align}
In terms of the Bezoutian this is equivalent to the
(non-symmetric) decomposition
\begin{equation}\label{eq:bezdec2}
B(a,b) = U_{\alpha_1}
\begin{pmatrix}
B(a_{-1},b_{-1}) & 0 \\
0 & 0
\end{pmatrix}
U_{\beta_1}^\top  + (\beta_1 - \alpha_1)
b_{\beta_1}a_{\alpha_1}^\top ,
\end{equation}
with $a_{\alpha_1}$ such that $a_{\alpha_1}^\top u_n(z)=a_{-1}(z)$ and $b_{\beta_1}$ likewise. \\
Iteration of this procedure gives
\begin{equation}\label{eq:algbezexp}
\frac{a(z)b(w)-a(w)b(z)}{z-w} = a(z)b(w) \sum_{k=1}^n (\beta_k-
\alpha_k)\frac{a_{-k}(w)b_{-k}(z)}{a_{-(k-1)}(z)b_{-(k-1)}(w)},
\end{equation}\label{eq:bezexp}
which gives the following expansion for the Bezout matrix
\begin{equation}\nonumber
B(a,b) = \sum_{k=1}^n (\beta_k-\alpha_k)U_{\alpha_1}\cdots
U_{\alpha_{k-1}}U_{\beta_{k+1}}\cdots U_{\beta_n}ee^\top
U_{\beta_1}^\top \cdots U_{\beta_{k-1}}^\top U^\top
_{\alpha_{k+1}}\cdots U^\top _{\alpha_n}.
\end{equation}
\end{proposition}
{\bf Proof.} Equation~(\ref{eq:bezdec1}) follows from elementary
computations. To prove~(\ref{eq:bezdec2}) we premultiply both
sides of the equation by $u_n(z)^\top$ and postmultiply them by
$u_n(w)$. The obtained left hand side then is obviously equal to
the left hand side of~(\ref{eq:bezdec1}). To show that the right
hand sides coincide one uses that $u_n(z)^\top
U_{\alpha_1}=(1-\alpha_1z)(u_{n-1}(z)^\top,0)+(0,\ldots,0,z^{n-1})$.
Then the assertion easily follows from the definition of
$B(a_{-1},b_{-1})$. To prove the other assertions, we proceed as
follows. First we show how the right hand sides of
equations~(\ref{eq:algbezexp}) and~(\ref{eq:bezexp}) are related.
We pre-multiply the right hand side of~(\ref{eq:bezexp}) by
$u_n(z)^\top$. The important key relation is
\[
u_n(z)^\top U_{\alpha_1}\cdots
U_{\alpha_{k-1}}U_{\beta_{k+1}}\cdots U_{\beta_n}e =
\prod_{j=1}^{k-1} (1-\alpha_jz)\prod_{j=k+1}^n (1-\beta_jz),
\]
which is easily shown to be true. Of course the right hand side of
this equation is nothing else but
\[
\frac{a(z)}{a_{-(k-1)}(z)}b_{-k}(z).
\]
Then post-multiplication of the obtained expression by $u_n(w)$
obviously results in the right hand side
of~(\ref{eq:algbezexp}).\\
We now show by induction that this is equal to $u_n(z)^\top
B(a,b)u_n(w)$. Let $A(z)=(1-\alpha_0z)a(z)$,
$B(z)=(1-\beta_0z)b(z)$. Define
\[ A_{-k}(z)=\frac{A(z)}{\prod_{j=0}^k(1-\alpha_jz)} \] and define $B_{-k}(z)$
likewise ($k=0,\ldots,n$). We also let $A_1(z)=A(z)$ and
$B_1(z)=B(z)$. We will use the following trivial identities. For
$k=1,\ldots,n$ we have $A_{-k}(z)=a_{-k}(z)$ and
$B_{-k}(z)=b_{-k}(z)$.
\begin{eqnarray*}
\lefteqn{A(z)B(w) \sum_{k=0}^n (\beta_k-
\alpha_k)\frac{A_{-k}(w)B_{-k}(z)}{A_{-(k-1)}(z)B_{-(k-1)}(w)} = }
\\
& &  (\beta_0-\alpha_0) a(w)b(z) + \\
& & (1-\alpha_0z)(1-\beta_0w) a(z)b(w) \sum_{k=1}^n (\beta_k-
\alpha_k)\frac{a_{-k}(w)b_{-k}(z)}{a_{-(k-1)}(z)b_{-(k-1)}(w)}.
\end{eqnarray*}
In view of~(\ref{eq:bezdec1}) and the induction assumption, this
equals  $\frac{A(z)B(w)-A(w)B(z)}{z-w}$. The proposition has been
proved. \hfill$\square$

\begin{corollary}\label{cor:bez0}
Let $\phi$ be a common zero of $\hat{a}$ and $\hat{b}$. Then
$a(z)=(1-\phi z)a_{-1}(z)$ and $b(z)=(1-\phi z)b_{-1}(z)$ and
\begin{equation}\label{eq:bez0}
B(a,b)=U_\phi
\begin{pmatrix}
B(a_{-1},b_{-1}) & 0 \\
0 & 0 \end{pmatrix} U_\phi^\top.
\end{equation}
\end{corollary}
{\bf Proof.} This is a straightforward consequence of the previous
proposition. \hfill $\square$\medskip\\
Proposition~\ref{prop:bezexp}  can  be used to show the well known
fact (see~\cite[theorem 13.1]{lancastertismenetsky}
or~\cite[theorem 8.4.3]{fuhrmann}) that the Bezout matrix $B(a,b)$
is non-singular iff $a$ and $b$ have no common factors. We use
corollary~\ref{cor:bez0} to give a description of the kernel of
the Bezout matrix.
\begin{corollary}
Let $\gamma_1,\ldots,\gamma_m$ be all the common zeros of
$\hat{a}$ and $\hat{b}$, with multiplicities $n_1,\ldots,n_m$. Let
$\ell$ be the last basis vector of $\openR^n$ and put $v^j_k =
(T_{\gamma_k}^jJ^{j-1})^\top \ell$ for $k=1,\ldots,m$ and
$j=1,\ldots,n_k$. Then $\ker B(a,b)$ is the linear span of the
vectors $v_k^j$.
\end{corollary}
{\bf Proof.} First we have to show that the vectors $v_k^j$ are
independent. Explicit computation of these vectors show that,
after multiplication with $P$, they are columns of the confluent
Vandermonde matrix associated with all zeros of $\hat{a}$, from
which independence then follows. For $j=1$, it follows immediately
from corollary~\ref{cor:bez0} that the $v_k^1$ belong to the
kernel of the Bezout matrix. When $\phi_k$ is common zero with
multiplicity $j>1$ we can factor the matrix $B(a_{-1},b_{-1})$
in~(\ref{eq:bez0}) like $B(a,b)$, but with one dimension less.
However, one can then show that also
\[
B(a,b)=U_{\phi_k}^2
\begin{pmatrix}
B(a_{-2},b_{-2}) & 0 \\
0 & 0
\end{pmatrix}
 (U_{\phi_k}^\top)^2,
\]
where for instance the $0$-matrix in the lower right corner now
has size $2\times 2$. Continuation of this procedure leads to
\[
B(a,b)=U_{\phi_k}^{j} \begin{pmatrix}
B(a_{-j},b_{-j}) & 0 \\
0 & 0 \end{pmatrix} (U_{\phi_k}^\top)^{j},
\]
for $j=1,\ldots,n_k$. Since the last $j$ columns of $B(a,b)$ are
thus zero vectors, one immediately sees that $B(a,b)v^j_k=0$. The
proof is complete upon noticing that the kernel of the Bezout
matrix has dimension equal to $n_1+\cdots +n_m$ (cf~\cite[theorem
8.4.3]{fuhrmann}). \hfill$\square$\medskip
\begin{remark}
{\em For more applications of confluent Vandermonde matrices to
the analysis of stationary ARMA processes, we refer
to~\cite{kleinspreijsiam}. }
\end{remark}

\section{The Bezoutian and the ARMA
polynomials}\label{sec:armabezout}

\setcounter{equation}{0}

In this section we continue to study some properties related to
the Bezout matrix, which (aimed at applications in
section~\ref{section:fisher}) we express in terms of the ARMA
polynomials $a$ and $c$ that define the process $y$ of
equation~(\ref{eq:arma}). For a polynomial
$a(z)=\sum_{k=0}^na_kz^k$ the matrix $S(a)$ is given by
\[
S(a)=
\begin{pmatrix}
  a_{1} & a_{2} & \cdots & a_{n} \\
  a_{2} & \cdots & a_{n} & 0 \\
  \vdots & \bdots & \bdots & \vdots \\
  a_{n} & 0 & \cdots & 0
\end{pmatrix}
\]
and $S(\hat{a})$ is given by
\[
S(\hat{a})=
\begin{pmatrix}
  a_{n-1} & a_{n-2} & \cdots & a_{0} \\
  a_{n-2} & \cdots & a_{0} & 0 \\
  \vdots & \bdots & \bdots & \vdots \\
  a_{0} & 0 & \cdots & 0
\end{pmatrix}.
\]
As before we will work with polynomials $a$ whose constant term
$a_0=1$. Notice that $S(\hat{a})$ is connected to the reciprocal
polynomial $\hat{a}$, $\hat{a}(z)=\sum_{k=0}^na_{n-k}z^k$, as is
$S(a)$ to $a$. Let $P$ be the `anti-diagonal identity' matrix in
$\openR^{n\times n}$, so with elements $P_{ij}=\delta_{i,n+1-j}$.
On a {\em Toeplitz} matrix $M$ pre- and postmultiplication by $P$
results in the same as transposition: $PMP =M^\top $. We will use
this property mainly for the choice $M=J$, the shift matrix.
\medskip\\
We continue under the assumption that the polynomials $\hat{a}$
and $\hat{c}$ have common degree $n$. One of the possible
relations between the Sylvester matrix $R(c,-a)$ and the Bezoutian
$B(c,a)$ is given below.
\begin{proposition}\label{prop:sylbez}
The matrices $R(c,-a)$ and $B(c,a)$ satisfy
\begin{equation}\nonumber
\begin{pmatrix}
P & 0 \\
PS(\hat{a})P & PS(\hat{c})P \end{pmatrix} R(c,-a)=
\begin{pmatrix}
I & 0 \\
0 & B(c,a) \end{pmatrix} \begin{pmatrix}
PS(\hat{c})P & S(c) \\
0 & I \end{pmatrix}.
\end{equation}

\end{proposition}
{\bf Proof.} This relation is just a variant on equation (21) on
page 460 of~\cite {lancastertismenetsky} and can be proven
similarly.
\hfill$\square$\medskip\\
We will use the short hand notation
\begin{equation}\label{eq:M}
M(c,a) =
\begin{pmatrix}
P & 0 \\
PS(\hat{a})P & PS(\hat{c})P \end{pmatrix}
\end{equation}
and
\begin{equation}\label{eq:N}
N(c) = \begin{pmatrix}
PS(\hat{c})P & S(c) \\
0 & I \end{pmatrix}.
\end{equation}
Notice that both $M(c,a)$ and $N(c)$ are
nonsingular if $a_0\neq 0$ and $c_0\neq 0$ (which is our case,
since we always work with $a_0=c_0=1$).

\begin{theorem}\label{thm:mam} Let $F$, $G$, $M(c,a)$ and $N(c)$ be as
in equations~(\ref{eq:F}), (\ref{eq:G}), (\ref{eq:M}) and
(\ref{eq:N}). The following identities hold true.
\begin{equation}
M(c,a)GM(c,a)^{-1}=\begin{pmatrix}
P(J-ea^\top )P & 0 \\
(c-a)e^\top  & PJP-ce^\top \end{pmatrix} =: G_M \label{eq:am}
\end{equation}
and

\begin{equation}
N(c)FN(c)^{-1}= \begin{pmatrix}
P(J-ea^\top )P & 0 \\
ee^\top        & J-ec^\top \end{pmatrix} =: F_N. \label{eq:an}
\end{equation}
Moreover we have the relation
\begin{equation}
G_M \begin{pmatrix}
I & 0 \\
0 & B(c,a) \end{pmatrix} = \begin{pmatrix}
I & 0 \\
0 & B(c,a) \end{pmatrix} F_N. \label{eq:amn}
\end{equation}
\end{theorem}
Before giving the proof of this theorem we formulate  a few
technical lemmas that will be of use in this proof.

\begin{lemma}\label{lemma:psp}
The following two equalities hold true.
\begin{equation}\nonumber
PS(\hat{c})P(J-ec^\top ) = (PJP-ce^\top )PS(\hat{c})P  =
\begin{pmatrix}
0      & \ldots  & 0      & 1       & 0 \\
\vdots &         &        & c_1     & \vdots  \\
0      &         &        & \vdots  & \vdots  \\
1      & c_1     & \cdots & c_{n-2} & 0 \\
0      & \cdots  & \cdots & 0       & -c_n
\end{pmatrix}.
\end{equation}
\end{lemma}
{\bf Proof.} Compare to the analogous statement in
\cite{lancastertismenetsky}, page 455. \hfill$\square$
\begin{lemma}
Let $g(z)=a(z)c(z)=\sum_{k=0}^{2n} g_kz^k$ and $g=(g^1,g^2)$, with
$g^1=(g_1,\ldots,g_n)$ and $g^2=(g_{n+1},\ldots,g_{2n})$. Then the
following identities hold true.
\begin{equation}\nonumber
S(\hat{c})Pe = e.
\end{equation}
\begin{equation}\nonumber
a^\top S(\hat{c})P = (g^1-c)^\top
\end{equation}
\begin{equation}
S(c)Pa = g^2. \label{eq:g2}
\end{equation}
\end{lemma}
{\bf Proof.} This is a straightforward verification.
\hfill$\square$\medskip\\
Along with the matrices $S(c)$ and $S(\hat{c})$ we also use the
Hankel matrix $\tilde{S}(\hat{c}) \in \openR^{(n+1)\times (n+1)}$
defined by
\begin{equation}\nonumber
\tilde{S}(\hat{c}) =
\begin{pmatrix}
c_n    & \ldots & \ldots & c_1 & 1 \\
\vdots &        &        & 1   & 0 \\
\vdots &        &        &     & \vdots \\
c_1    & 1      &        &     & \vdots \\
1      & 0      & \ldots & \ldots & 0
\end{pmatrix}.
\end{equation}

\begin{lemma} \label{lemma:symm}
One has
\begin{equation}\nonumber
JS(\hat{c})+ec^\top P =
\begin{pmatrix}
I & 0 \end{pmatrix} \tilde{S}(\hat{c}) \begin{pmatrix}
I \\
0 \end{pmatrix}.
\end{equation}
In particular the matrix $JS(\hat{c})+ec^\top P$ is symmetric.
\end{lemma}
{\bf Proof.} The following relations are immediate.
\begin{equation}\label{eq:ss1}
S(\hat{c}) = \begin{pmatrix} 0 & I \end{pmatrix}
\tilde{S}(\hat{c}) \begin{pmatrix}
I \\
0 \end{pmatrix}
\end{equation}
and
\begin{equation}\label{eq:ss2}
ec^\top P = \begin{pmatrix} e & 0 \end{pmatrix} \tilde{S}(\hat{c})
\begin{pmatrix}
I \\
0 \end{pmatrix}.
\end{equation}
Use equations~(\ref{eq:ss1}) and~(\ref{eq:ss2}) to write
\begin{eqnarray*}
JS(\hat{c})+ec^\top P & = & J\begin{pmatrix} 0 & I \end{pmatrix}
\tilde{S}(\hat{c}) \begin{pmatrix}
I \\
0 \end{pmatrix} + \begin{pmatrix} e & 0 \end{pmatrix}
\tilde{S}(\hat{c}) \begin{pmatrix}
I \\
0
\end{pmatrix} \\
& = & \big(\begin{pmatrix} 0 & J \end{pmatrix} + \begin{pmatrix} e
& 0 \end{pmatrix})\big) \tilde{S}(\hat{c}) \begin{pmatrix}
I \\
0 \end{pmatrix},
\end{eqnarray*}
from which the result follows.
\hfill$\square$\medskip \\
{\bf Proof of theorem~\ref{thm:mam}.} Compute the two products
$M(c,a)G$ and $G_MM(c,a)$ to get respectively
\begin{equation}
\begin{pmatrix}
P(J-ea^\top )    & 0 \\
PS(\hat{a})P(J-ea^\top ) & PS(\hat{c})P(J-ec^\top ) \end{pmatrix}
\label{eq:first}
\end{equation}
and
\begin{equation}
\begin{pmatrix}
P(J-ea^\top )    & 0 \\
(c-a)e^\top  + (PJP-ce^\top )PS(\hat{a})P & (PJP-ce^\top
)PS(\hat{c})P \end{pmatrix}. \label{eq:second}
\end{equation}
Clearly we only have to look at the 21- and 22-blocks. Comparing the
22-blocks is just the content of lemma~\ref{lemma:psp} \\
We focus on the 21-blocks. Use lemma~\ref{lemma:psp} again to write
$PS(\hat{a})P(J-ea^\top ) $ (the 21-block of (\ref{eq:first}) as
\[
(PJP- ae^\top )PS(\hat{a})P = PJS(\hat{a})-ae^\top P =
PJS(\hat{a})-ce^\top P +(c-a)e^\top P = \]
\[ PJS(\hat{a})-ce^\top PS(\hat{a})P +(c-a)e^\top P, \]
which is just the 21-block of (\ref{eq:second}). This proves the
identity~(\ref{eq:am}). \\
Next we prove (\ref{eq:an}). Write $F = \begin{pmatrix}
J-eg^{1^\top}  & -eg^{2^\top}  \\
ee^\top P     & J \end{pmatrix}$. Work out the products $N(c)F$
and $F_NN(c)$ to get respectively
\begin{equation}
\begin{pmatrix}
PS(\hat{c})P(J-eg^{1^\top} )+S(c)ee^\top P & -PS(\hat{c})Peg^{2^\top} + S(c)J \\
ee^\top P                             & J \end{pmatrix}
\label{eq:first1}
\end{equation}
and
\begin{equation}
\begin{pmatrix}
P(J-ea^\top )S(\hat{c})P & P(J-ea^\top )PS(c) \\
ee^\top PS(\hat{c})P     & ee^\top S(c)+J-ec^\top \end{pmatrix}.
\label{eq:second1}
\end{equation}
Compare now the corresponding blocks in these two matrices. We start with the
11-block of (\ref{eq:first1}).
Write it as
\[
P(S(\hat{c})PJP + Pce^\top )P - Peg^{1^\top}  = P(S(\hat{c})J^\top
+ Pce^\top )P - Peg^{1^\top}
\]
and use the symmetry asserted in lemma~\ref{lemma:symm} to get
\[ P(JS(\hat{c}) + ec^\top P)P -
Peg^{1^\top}  =  PJS(\hat{c})P + e(c^\top  -g^{1^\top} ) =
PJS(\hat{c})P - Pea^\top S(\hat{c})P,
\]
which equals the 11-block of (\ref{eq:second1}). \\
Next we consider the 12-blocks. Start with (\ref{eq:first1}):
\[
-PS(\hat{c})Peg^{2^\top}  + S(c)J =-Peg^{2^\top}  + S(c)J =
-Pea^\top PS(c) +S(c)J,
\]
where the last equality just follows from (\ref{eq:g2}).
Since $S(c)J$ is symmetric it is equal to
$J^\top S(c)=PJPS(c)$. Hence
\[
-Pea^\top PS(c) +S(c)J=P(J-ea^\top )PS(c),
\] which is equal
to the 12-block of (\ref{eq:second1}).
Comparison of the other blocks is trivial. \\
Finally we prove (\ref{eq:amn}). Remember that $GR(c,-a)=R(c,-a)F$
(proposition~\ref{prop:eps2}). Write $\bar{B}$ for
$\begin{pmatrix}
I & 0 \\
0 & B(c,a) \end{pmatrix} $. Then we have the string of equalities
\[ G_M\bar{B}=G_MM(c,a)R(c,-a)N_c^{-1} = M(c,a)GR(c,-
a)N_c^{-1} = \]
\[ M(c,a)R(c,-a)FN_c^{-1} = \bar{B}N_cFN_c^{-1} =
\bar{B}F_N.
\]
This proves the last assertion of the theorem. \hfill $\square$

\section{Stein equations}\label{sec:stein}
\setcounter{equation}{0}

We start this section with considering two Stein equations that
involve the matrices $F_M$ and $G_N$ of equations~(\ref{eq:am})
and~(\ref{eq:an}).
\begin{proposition}\label{prop:HQ}
Let $e_P^\top =[e^\top P,0]^\top $ and let $H(\theta)$ and
$Q(\theta)$ be the unique solutions  to the following Stein
equations
\begin{eqnarray}
H & = & G_MHG_M^\top  + e_Pe_P^\top  \label{eq:H}\\
Q & = & F_NQF_N^\top  + e_Pe_P^\top .\label{eq:Q}
\end{eqnarray}
Then $Q(\theta)$ is strictly positive definite. Moreover,
$H(\theta)$ and $Q(\theta)$ are related by
\begin{equation}\label{eq:GQ}
H(\theta) = \begin{pmatrix}
I & 0 \\
0 & B(c,a) \end{pmatrix} Q(\theta) \begin{pmatrix}
I & 0 \\
0 & B(c,a) \end{pmatrix}.
\end{equation}
\end{proposition}
{\bf Proof.} To show that $Q(\theta)$ is strictly positive
definite, it is sufficient to show that the pair $(F_N,e_P)$ is
controllable. Let $T=\begin{pmatrix}P & 0 \\ 0 & I
\end{pmatrix}$. For computational reasons it is more convenient to
show controllability of the pair $(A,b)$, where $A= TF_NT^{-1}$
and $b=TF_N$. Observe that $b$ is the first standard basis vector
in $\openR^{2n}$, whereas $A=\begin{pmatrix}
J-ea^\top P & 0 \\
ee^\top P       & J-ec^\top
\end{pmatrix}$. If one computes the controllability matrix
$(b, Ab,\ldots, A^{2n-1}b)$, then standard calculations lead to
the fact that this matrix is upper triangular with only ones on
the diagonal. Hence it has full rank. By theorem 8d.66
of~\cite{callierdesoer}
the matrix $Q(\theta)$ is strictly positive definite. \\
Multiply equation~(\ref{eq:Q}) with $Q=Q(\theta)$ on the right and
on the left by the symmetric matrix $T=\begin{pmatrix} I & 0 \\ 0
& B(c,a)
\end{pmatrix}$ and put $H=TQ(\theta)T$. In view of
relation~(\ref{eq:amn}) we then obtain equation~(\ref{eq:H}).
Hence $H$ must be equal to $H(\theta)$ by uniqueness of the
solution. This shows~(\ref{eq:GQ}). \hfill$\square$

\begin{corollary}\label{cor:Gac}
The matrix $H(\theta)$ is non-singular iff the polynomials $a$ and
$c$ have no common factors.
\end{corollary}
{\bf Proof.} The matrix $Q(\theta)$ is non-singular and $B(c,a)$
is singular iff the polynomials $a$ and $c$ have no common
factors. \hfill$\square$
\begin{remark}
{\em  If the polynomials $a$ and $c$ have a common factor $(1-\phi
z)$, then the expression for $B(c,a)$ of equation~(\ref{eq:bez0})
can be applied to obtain a rank factorization of $H(\theta)$. }
\end{remark}
Along with the matrices $H(\theta)$ and $Q(\theta)$ of
proposition~\ref{prop:HQ} we will also work with the matrices
$I(\theta)$ and $P(\theta)$, defined in
\begin{align}
I(\theta) & = M(c,a)^{-1}H(\theta)M(c,a)^{-\top} \label{eq:IG} \\
P(\theta) & = N(c)^{-\top}Q(\theta)N(c)^{-1}, \label{eq:PQ}
\end{align}
where $M(c,a)$ and $N(c)$ are as in~(\ref{eq:M}) and~(\ref{eq:N}).
In view of equations~(\ref{eq:am}) and~(\ref{eq:an}) and
proposition~\ref{prop:HQ}, we have that $I(\theta)$ and
$P(\theta)$ are solutions to the Stein equations
\begin{align}\label{eq:I}
I & = FIF^\top +  BB^\top \\
P & = GPG^\top + ee^\top.
\end{align}
\begin{corollary}\label{cor:PI}
The matrix $P(\theta)$ is non-singular and the matrix $I(\theta)$
is non-singular iff the polynomials $a$ and $c$ have no common
factors.
\end{corollary}
{\bf Proof.} Since the matrices $M(c,a)$ and $N(c)$ are
non-singular, the results follows from proposition~\ref{prop:HQ}
and corollary~\ref{cor:Gac}.\hfill$\square$

\section{Fisher's information matrix}\label{section:fisher}

\setcounter{equation}{0}

We consider again the ARMA process $y$ defined by~(\ref{eq:arma}).
Let the variance of the white noise sequence be $\sigma^2$. Assume
that the process is stationary. Traditionally, the Fisher
information matrix is defined as the covariance matrix of the
score function. Let us assume that the process is also Gaussian
and that $\sigma^2$ is a known constant. If we have observations
$y_1,\ldots,y_n$, the log likelihood $\ell_n(\theta)$ is then
essentially given by
\[
\ell_n(\theta)=\frac{n}{2}\log \sigma^2
-\frac{1}{2\sigma^2}\sum_{t=1}^n\eps_t(\theta)^2.
\]
The score function, by definition the gradient of
$\ell_n(\theta)$, is then given by
\begin{equation}\label{eq:score}
\dot{\ell}_n(\theta)=-\frac{1}{\sigma^2}\sum_{t=1}^n\eps_t(\theta)\dot{\eps}_t(\theta).
\end{equation}
Remember that for $\dot{\eps}_t(\theta)$ we have the
realizations~(\ref{eq:obsreal}) and~(\ref{eq:epsstate}). In
remark~\ref{rem:independent} we observed that
$\dot{\eps}_t(\theta)$ is (stochastically) independent of
$\eps_t$, and we also have $\E \dot{\ell}_t(\theta)=0$
(expectation taken w.r.t.~the distribution under $\theta$). The
covariance matrix $I_n(\theta,\sigma^2)$, by definition the
covariance matrix of $\dot{\ell}_n(\theta)$, can then be computed
as
\[
I_n(\theta,\sigma^2)=\cov(\dot{\ell}_n(\theta)) =\E
\dot{\ell}_n(\theta)^\top\dot{\ell}_n(\theta).
\]
It then follows from~(\ref{eq:score}) and the above mentioned
independence that
\[
I_n(\theta,\sigma^2)=\frac{1}{\sigma^4}\sum_{t=1}^n \E
\eps_t(\theta)^2
\E\dot{\eps}_t(\theta)^\top\dot{\eps}_t(\theta)=\frac{1}{\sigma^2}\sum_{t=1}^n
\E\dot{\eps}_t(\theta)^\top\dot{\eps}_t(\theta).
\]
For
$i_t(\theta,\sigma^2):=\E\dot{\eps}_t(\theta)^\top\dot{\eps}_t(\theta)$
we get from equation~(\ref{eq:epsstate}) and the independence of
$\eps_t(\theta)$ and $\dot{\eps}_t(\theta)$ the recursion
\begin{equation}\label{eq:it}
i_{t+1}(\theta,\sigma^2)=Fi_t(\theta,\sigma^2)F^\top+\sigma^2BB^\top.
\end{equation}
Under the stationarity assumption we have
$i_{t+1}(\theta,\sigma^2)=i_{t}(\theta,\sigma^2)$ and we simply
write $i(\theta,\sigma^2)$. Hence
$I_n(\theta,\sigma^2)=nI(\theta,\sigma^2)$, where
\begin{equation}\label{eq:its}
I(\theta,\sigma^2)=\frac{1}{\sigma^2}\E i(\theta,\sigma^2).
\end{equation}
Without stationary initial conditions, but still with $\eps$ a
Gaussian white noise process, one can show (but this is not
relevant for the present paper) that
$\frac{1}{n}I_n(\theta,\sigma^2)\to I(\theta,\sigma^2)$. Hence the
matrix $I(\theta,\sigma^2)$ is also relevant in a non-stationary
situation. We call $I(\theta,\sigma^2)$ the asymptotic Fisher
information matrix. Summing up intermediate results, we obtain the
following theorem.
\begin{theorem}\label{thm:fsing}
The asymptotic Fisher information matrix $I(\theta,\sigma^2)$ of
the ARMA process defined by~(\ref{eq:arma})  is the same as the
matrix $I(\theta)$ defined in equation~(\ref{eq:IG}). Hence it is
the unique solution to the Stein equation
\begin{equation}\label{eq:II}
I=FIF^\top+BB^\top,
\end{equation}
and thus independent of $\sigma^2$. Moreover this matrix is
non-singular iff the polynomials $a$ and $c$ have no common
factors.
\end{theorem}
{\bf Proof.} From equations~(\ref{eq:it}) and~(\ref{eq:its}) and
the stationarity assumption, one immediately sees that
$I(\theta,\sigma^2)$ satisfies~(\ref{eq:II}), which is just
equation~(\ref{eq:I}). Hence the matrices $I(\theta,\sigma^2)$ and
$I(\theta)$ are the same and the characterization of
non-singularity is nothing else but corollary~\ref{cor:PI}.
\hfill$\square$\medskip\\
The conclusion of this theorem has been proved in~\cite{gyor} by
different means, involving representations of the Fisher
information matrix as an integral in the complex plane and the
following lemma of which we give an alternative proof.
\begin{lemma}\label{lemma:fact}
Let $I(\theta)$ be the Fisher information matrix and $P(\theta)$
as in~(\ref{eq:PQ}). Then the following factorization holds.
\begin{equation}
I(\theta )=R(c,-a)P(\theta )R(c,-a)^\top   \label{eq:fisherp}
\end{equation}
\end{lemma}
{\bf Proof.} This follows from proposition~\ref{prop:sylbez},
combined with equations~(\ref{eq:GQ}), (\ref{eq:IG})
and~(\ref{eq:PQ}).
\hfill$\square$\medskip\\
Since the matrix $P(\theta)$ is non-singular, also
lemma~\ref{lemma:fact} illustrates the fact that $I(\theta)$ is
non-singular iff $a$ and $c$ have no common factors. Moreover,
looking at equation~(\ref{eq:I}), we see that $I=I(\theta)$ is
non-singular iff the pair $(F,B)$ is controllable. But the
controllability matrices $\mathcal{R}(F,B)$ and $\mathcal{R}(G,e)$
satisfy the easily verified relation
$\mathcal{R}(F,B)=R(c,-a)\mathcal{R}(G,e)$. Since the matrix
$\mathcal{R}(G,e)$ has full rank, we see that $(F,e)$ is
controllable iff $R(c,-a)$ is non-singular, which leads to another
way of showing the conclusion of theorem~\ref{thm:fsing}.


\begin{thebibliography}{99}
\bibitem{astromsoderstrom}
K.J.~\AA str\"{o}m and T. S\"{o}derstr\"{o}m (1974), Uniqueness of
the Maximum Likelihood Estimates of the Parameters of an ARMA
model, {\em IEEE Transactions on Automatic Control} {\bf 19},
769--773.

\bibitem{barnett}
S.~Barnett (1971), A new formulation of the theorems of Hurwitz,
Routh and Sturm, {\em J.\ Inst.\ Maths.\ Applics}. {\bf 8},
240--250.

\bibitem{barnett2}   S.\ Barnett (1973), Matrices, polynomials
and linear time invariant systems, {\em IEEE Trans.\ Automat.\
Control} {\bf 18}, 1--10.

\bibitem{callierdesoer}
F.M.~Callier and C.A.~Desoer (1991), {\em Linear System Theory},
Springer.

\bibitem{cramer}   H.~Cram\'{e}r (1951), {\em Mathematical
methods of statistics}, Princeton University Press.

\bibitem{fuhrmann}
P.A.~Fuhrmann (1996), {\em A Polynomial Approach to Linear
Algebra}, Springer.

\bibitem{kalman}   R.E.~Kalman (1963), Mathematical description
of linear dynamical systems, {\em SIAM J.\ Contr}. {\bf 1},
152--192.


\bibitem{kleinspreij}   A.~Klein and P.J.C.~Spreij (1996), On
Fisher's Information Matrix of an ARMAX Process and Sylvester's
Resultant Matrices, {\em Linear Algebra Appl.} {\bf 237/238},
579--590.


\bibitem{gyor}
A.~Klein and P.J.C.~Spreij (1997), On Fisher's information matrix
of an ARMA process, {\em Stochastic differential  and difference
equations}, I.~Csiszar and Gy.~Michaletzky eds., 273--284,
Birkh\"auser.

\bibitem{kleinspreijsiam}
A.~Klein and P.J.C.~Spreij (2003), Some results on Vandermonde
matrices with an application to time series analysis, {\em SIAM
Journal on Matrix Analysis} {\bf 25(1)}, 213--223.



\bibitem{lancastertismenetsky}   P.~Lancaster and M.~Tismenetsky (1985),
{\em The Theory of Matrices}, Academic Press.


\bibitem{mcleod}   A.I.~McLeod (1993), A note on ARMA model
redundancy, {\em J. of Time Series Analysis} {\bf 14(2)},
207--208.

\bibitem{potscher}   B.M.~P\"{o}tscher (1985), The behaviour of
the Lagrangian multiplier test in testing the orders of an ARMA model, {\em %
Metrika} {\bf 32},  129--150.



\bibitem{rao}   C.R.~Rao (1965), {\em Linear statistical
inference and its applications},  Wiley.



\bibitem{rothenberg}   T.J.~Rothenberg (1971), Identification in
parametric models, {\em Econometrica} {\bf 39}, 577--591.

\bibitem{soderstrom}   T.~S\"{o}derstr\"{o}m and P.\ Stoica
(1989), {\em System identification}, Prentice hall.

\bibitem{vanderwaerden}   B.L.~van der Waerden (1966), {\em %
Algebra I}, Springer.
\end{thebibliography}
\end{document}